\documentclass[12pt]{article}

\usepackage{amsfonts}

\newtheorem{thm}{Theorem}[section]
\newtheorem{lemma}[thm]{Lemma}
\newtheorem{cor}[thm]{Corollary}

\newtheorem{prop}[thm]{Proposition}
\newtheorem{conjecture}{Conjecture}

\newcommand{\proof
}{\par\medskip\noindent {\bf Proof.\ \ }}

\newcommand{\be}{\begin{equation}}
\newcommand{\ee}{\end{equation}}
\newcommand{\openbox}{\leavevmode
  \hbox to8pt{\hfil\vrule\vbox to6pt{\hrule width6pt\vfil\hrule}\vrule}}

\newcommand{\qed}{\hbox to5pt{ } \hfill \openbox\bigskip\medskip}

\newcommand{\cF}{\mbox{$\cal F$}}
\newcommand{\cG}{\mbox{$\cal G$}}

\title{Sunflowers and $L$-intersecting families}
\author{G\'abor Heged\H{u}s
\\{\normalsize  \'Obuda University}
}
\begin{document}
\maketitle
\begin{abstract}
Let $f(k,r,s)$ stand for the least number so that if $\cF$ is an arbitrary $k$-uniform, $L$-intersecting set system, where $|L|=s$, and  $\cF$ has more than  $f(k,r,s)$ elements, then   $\cF$ contains a sunflower with $r$ petals.

We give an upper bound for $f(k,3,s)$.

Let $g(k,r,\ell)$ be the least number so that any $k$-uniform, $\ell$-intersecting set system of  more than  $g(k,r,\ell)$ sets contains a sunflower with $r$ petals. 

We give also an upper bound for $g(k,r,\ell)$.
\end{abstract}


\medskip
\section{Introduction}

Let  $[n]$ stand for the set $\{1,2,
\ldots, n\}$. We denote the family of all subsets of $[n]$  by $2^{[n]}$. 

Let $X$ be a fixed subset of $[n]$. For an integer $0\leq k\leq n$ we denote by
${X \choose k}$ the family of all  $k$ element subsets of $X$.

We call a family $\cF$ of subsets of $[n]$ {\em $k$-uniform}, if $|F|=k$ for each $F\in \cF$.
 
A family $\cF=\{F_1,\ldots ,F_m\}$ of subsets of $[n]$ is a {\em sunflower} (or {\em $\Delta$-system}) with $m$ petals if
$$
F_i\cap F_j=\bigcap\limits_{t=1}^m F_t
$$
for each $1\leq i,j\leq m$.

The intersection of the members of a sunflower form its {\em kernel}.
Clearly a a family of disjoint sets is a sunflower with empty kernel.

Erd\H{o}s and  Rado gave an upper bound for the size of a $k$-uniform family without a sunflower with $r$ petals in \cite{ER}.
\begin{thm} \label{Sthm} (Sunflower theorem)
If $\cF$ is a $k$-uniform set system with more than 
$$
k!(r-1)^k\big( 1-\sum_{t=1}^{k-1} \frac{t}{(t+1)! (r-1)^t}\big)
$$ 
members, then $\cF$ contains a sunflower with $r$ petals.
\end{thm}
                  
Kostochka improved this upper bound in \cite{K}. 
\begin{thm} \label{Kost}
Let $r>2$ and $\alpha>1$ be fixed integers. Let $k$ be an arbitrary integer. Then there exists a constant $D(r,\alpha)$ such that if  $\cF$ is a $k$-uniform set system with more than 
$$
D(r,\alpha)k! \Big( \frac{(\log \log \log k)^2}{\alpha \log \log k}\Big)^k
$$
 members, then $\cF$ contains a sunflower with $r$ petals.
\end{thm}

Erd\H{o}s and  Rado gave in \cite{ER} a construction of a $k$-uniform  set system with $(r-1)^k$ members such that $\cF$ does not contain any sunflower with $r$ petals. Later Abbott, Hanson and Sauer improved this construction in \cite{AHS} and proved the following result.

\begin{thm} \label{AHS}
There exists a $c>0$ positive constant and a $k$-uniform set system $\cF$ such that 
$$
|\cF|> 2\cdot 10^{k/2-c\log k}
$$
and $\cF$ does not contain any sunflower with $3$ petals.
\end{thm}

Erd\H{o}s and  Rado conjectured also the following statement in \cite{ER}.
\begin{conjecture} \label{Econj}
For each $r$, there exists a constant $C_r$ such that  if  $\cF$ is a $k$-uniform set system with more than 
$$
C_r^k
$$
members,  then $\cF$ contains a sunflower with $r$ petals.
\end{conjecture}
Erd\H{o}s has offered 1000 dollars for the proof or disproof of this conjecture for $r=3$ (see \cite{E}).

We prove here  Conjecture \ref{Econj} in the case of some special $L$-intersecting and $\ell$-intersecting families.

A family $\cF$ is {\em $\ell$-intersecting}, if $|F\cap F'|\geq \ell$ whenever $F, F'\in \cF$. Specially, $\cF$ is an {\em intersecting} family, if $F\cap F'\ne \emptyset$ whenever $F, F'\in \cF$.

Erd\H{o}s, Ko and  Rado proved the following well-known result in \cite{EKR}:
\begin{thm} \label{EKZ}
Let $n,k,t$ be integers with $0<t<k< n$ . Suppose $\cF$ is a  $t$-intersecting,  $k$-uniform family of subsets of $[n]$. Then for $n>n_0(k,t)$,
$$
|\cF|\leq {n-t \choose k-t}.
$$
Further,  $|\cF|= {n-t \choose k-t}$ if and only if for some $T\in {[n]\choose t}$ we have 
$$
\cF=\{F\in {[n]\choose k}:~ T\subseteq F\}.
$$ 
\end{thm}

Let $L$ be a set of nonnegative integers. A family $\cF$ is {\em $L$-intersecting}, if $|E\cap F|\in L$ for every pair $E,F$ of distinct members of $\cF$. In this terminology a $k$-uniform $\cF$ set system is a $t$-intersecting family iff  it is an $L$-intersecting family, where $L=\{t, t+1, \ldots , k-1\}$.

The following result gives a remarkable upper bound for the size of a $k$-uniform $L$-intersecting family (see \cite{RW}).

\begin{thm} \label{RW} (Ray-Chaudhuri--Wilson) Let $0<s\leq k\leq n$ be positive integers. Let $L$ be a set of $s$ nonnegative integers and $\cF$  an $L$-intersecting, $k$-uniform family  of subsets of $[n]$. Then
$$
|\cF|\leq  {n \choose s}.
$$
\end{thm}

Deza proved the following result in \cite{D}.
\begin{thm} (Deza) \label{Deza}
Let $\lambda>0$ be a positive integer. Let $L:=\{\lambda\}$. If $\cF$ is an $L$-intersecting, $k$-uniform family  of subsets of $[n]$, then either
$$
|\cF|\leq k^2-k+1
$$
or $\cF$ is a sunflower, i.e. all the pairwise intersections are the same set with $\lambda$ elements.
\end{thm}

We generalize Theorem \ref{Deza} for $L$-intersecting families in the following.

\begin{thm} \label{main}
Let $\cF$ be a  family  of subsets of $[n]$ such that $\cF$ does not contain any sunflowers with three petals. Let $L=\{\ell_1 <\ldots <\ell_s\}$ be a set of $s$ non-negative integers. Suppose that $\cF$ is a $k$-uniform, $L$-intersecting family. Then
$$
|\cF|\leq (k^2-k+2)8^{(s-1)}2^{(1+\frac{\sqrt{5}}{5})k(s-1)}.
$$
\end{thm}

Next we improve Theorem \ref{Sthm} in the case of $\ell$-intersecting families. 

\begin{thm} \label{main2}
Let $r>2$ and $\alpha>1$ be fixed integers. Let $\cF$ be an  $\ell$-intersecting, $k$-uniform family  of subsets of $[n]$ such that $\cF$ does not contain any sunflowers with $r$ petals. Then there exists a constant $D(r,\alpha)$ such that 
$$
|\cF|\leq D(r,\alpha){k\choose \ell} (k-\ell)! \Big( \frac{(\log \log \log (k-\ell))^2}{\alpha \log \log (k-\ell)}\Big)^{k-\ell}.
$$
\end{thm}

\begin{cor} \label{main3}
Let $r>2$, $k>1$ be fixed integers. Let $\ell:= \lceil k-\frac{k}{\log k}\rceil$. Let $\cF$ be an $\ell$-intersecting,  $k$-uniform family  of subsets of $[n]$ such that $\cF$ does not contain any sunflowers with $r$ petals. Then there exists a constant $D(r)$ such that
$$
|\cF|\leq D(r)4^k.
$$
\end{cor}

We present our proofs in Section 2. We give some concluding remarks in Section 3. 

\section{Proofs}

We start our proof with an elementary fact.

\begin{lemma} \label{help}
Let $0\leq r\leq n$ be integers. Then 
$$
{n \choose r}\leq {n-1 \choose r+1}
$$
if and only if
$$
r^2+(1-3n)r+n^2-2n \geq 0.
$$
\end{lemma}
\qed

\begin{cor} \label{help2}               
Let $0\leq r\leq n$ be integers. If $0\leq r \leq\frac{3n-1-\sqrt{5}(n+1)}{2}$, then 
$$
{n \choose r}\leq {n-1 \choose r+1}
$$
\end{cor}
\qed

We use the following easy Lemma in the  proof of our main results.

\begin{lemma} \label{lem_binom}
Let $0\leq \ell \leq k-1$ be integers. Then
$$
{2k-\ell \choose \ell +1}\leq 8\cdot 2^{(1+\frac{\sqrt{5}}{5})k}.
$$
\end{lemma}
\proof
First suppose that
$$
\ell \leq \lceil(1- \frac{\sqrt{5}}{5})k-(1+\frac{2\sqrt{5}}{5}) \rceil.
$$
Then
$$
{2k-\ell \choose \ell +1}\leq {2k- \lceil (1- \frac{\sqrt{5}}{5})k-(1+\frac{2\sqrt{5}}{5}) \rceil \choose \lceil (2- \frac{\sqrt{5}}{5})k-(1+\frac{2\sqrt{5}}{5}) \rceil } \leq 
$$
$$
\leq 2^{2k- \lceil (1- \frac{\sqrt{5}}{5})k-(1+\frac{2\sqrt{5}}{5}) \rceil }\leq 8\cdot 2^{(1+\frac{\sqrt{5}}{5})k}.
$$
The first inequality follows easily from Corollary \ref{help2}. Namely if 
$$
\ell \leq \lceil(1- \frac{\sqrt{5}}{5})k-(1+\frac{2\sqrt{5}}{5}) \rceil, 
$$
then 
$$
\ell +1\leq \lceil \frac{3-\sqrt{5}}{2}(2k-\ell)-\frac{1+\sqrt{5}}{2}   \rceil
$$
and we can apply Corollary \ref{help2} with the choices $r:=\ell +1$ and $n:=2k-\ell$.

Secondly,  suppose that
$$
\ell>  \lceil(1- \frac{\sqrt{5}}{5})k-(1+\frac{2\sqrt{5}}{5}) \rceil.
$$
Then
$$
2k-\ell \leq 2k-  \lceil(1- \frac{\sqrt{5}}{5})k-(1+\frac{2\sqrt{5}}{5}) \rceil \leq \lceil 2+(1+\frac{\sqrt{5}}{5})k \rceil, 
$$
hence
$$
{2k-\ell \choose \ell +1}\leq {2k- \lceil (1- \frac{\sqrt{5}}{5})k-(1+\frac{2\sqrt{5}}{5}) \rceil \choose \ell +1} \leq 
$$
$$
\leq 2^{2k- \lceil (1- \frac{\sqrt{5}}{5})k-(1+\frac{2\sqrt{5}}{5}) \rceil } \leq 8\cdot 2^{(1+\frac{\sqrt{5}}{5})k}.
$$

\qed

The soul of the proof of our main result is the following Lemma.

\begin{lemma} \label{soul}
Let $\cF$ be an  $\ell$-intersecting, $k$-uniform family  of subsets of $[n]$ such that $\cF$ does not contain any sunflowers with three petals. Suppose that there exist $F_1, F_2\in \cF$ distinct subsets such that $|F_1\cap F_2|=\ell$. Let $M:=F_1\cup F_2$. Then 
$$
|F\cap M|>\ell
$$
for each $F\in \cF$. 
\end{lemma}
\proof

Clearly $F\cap F_1 \subseteq F\cap M$ for each $F\in \cF$, hence 
$$
|F\cap M|\geq \ell
$$
for each $F\in \cF$. 

We prove by an indirect argument. Suppose that there exists an $F\in \cF$ such that $|F\cap M|=\ell$. Clearly $F\ne F_1$ and $F\ne F_2$. Let $G:=F_1 \cap F_2$. Then $|G|=\ell$ by assumption. It follows from $F\cap F_1 \subseteq F\cap M$ and 
$\ell\leq |F\cap F_1|\leq |F\cap M|=\ell$  that $F\cap F_1=F\cap M$. Similarly  $F\cap F_2=F\cap M$. Consequently $F\cap F_1=F\cap F_2$. We get that $F\cap F_2=F\cap G=F\cap F_1$. Since 
$\ell=|F\cap F_1|=|F\cap G|\leq |G|=\ell$ and $F\cap G\subseteq G$, hence $G=F\cap G=F\cap F_2=F\cap F_1$, so $\{F,F_1,F_2\}$ is a sunflower with three petals, a contradiction. \qed

{\bf Proof of Theorem \ref{main}:}

We apply  induction on $|L|=s$. If $s=1$, then our result follows from Theorem \ref{Deza}.

Suppose that Theorem \ref{main} is true for $s-1$ and now we attack the case $|L|=s$.

If $|F\cap F'|\ne \ell_1$ holds for each distinct $F, F'\in \cF$, then $\cF$ is actually an $L':=\{\ell_2, \ldots ,\ell_s\}$-intersecting system and the much stronger upper bound 
$$
|\cF|\leq (k^2-k+2)8^{(s-2)}2^{(1+\frac{\sqrt{5}}{5})k(s-2)}.
$$ 
follows from the induction.

Hence we can suppose that there exist $F_1, F_2\in \cF$ such that $|F_1\cap F_2|= \ell_1$. Let $M:=F_1 \cup F_2$. Clearly $\cF$ is an $\ell_1$-intersecting family. It follows from Lemma \ref{soul} that 
\begin{equation} \label{inter}
|F\cap M|>\ell_1
\end{equation}
for each $F\in \cF$. Clearly $|M|=2k-\ell_1$.

Let $T$ be a fixed subset of $M$ such that $|T|=\ell_1+1$.
Define the family
$$
\cF(T):=\{F\in \cF:~ T\subseteq M\cap F\}. 
$$
Let $L':=\{\ell_2, \ldots, \ell_s\}$. Clearly $|L'|=s-1$.
Then $\cF(T)$ is an $L'$-intersecting, $k$-uniform family, because $\cF$ is an $L$-intersecting family and $T\subseteq F$ for each $F\in\cF(T)$. The following Proposition follows easily from (\ref{inter}).
\begin{prop} \label{Funion}
$$
\cF=\bigcup\limits_{T\subseteq M, |T|=\ell_1+1} \cF(T).
$$
\end{prop}

\qed

Let $T$ be a fixed, but arbitrary subset of $M$ such that $|T|=\ell_1+1$. Consider the set system
$$
\cG(T):=\{F\setminus T:~ F\in \cF(T)\}.
$$
Clearly $|\cG(T)|=|\cF(T)|$. Let $\overline{L}:=\{\ell_2-\ell_1-1, \ldots, \ell_s-\ell_1-1\}$. Here $|\overline{L}|=s-1$.
Since  $\cF(T)$ is an $L'$-intersecting, $k$-uniform family, thus $\cG(T)$ is an $\overline{L}$-intersecting, $(k-\ell_1-1)$-uniform family. It follows from the inductional hypothesis that
$$
|\cF(T)|=|\cG(T)|\leq (k^2-k+2)8^{(s-2)}2^{(1+\frac{\sqrt{5}}{5})k(s-2)}.
$$

Finally Proposition \ref{Funion} implies that
$$
|\cF|\leq \sum_{T\subseteq M, |T|=\ell_1+1} |\cF(T)|\leq {2k-\ell_1 \choose \ell_1+1} (k^2-k+2)8^{(s-2)}2^{(1+\frac{\sqrt{5}}{5})k(s-2)}.
$$
But 
$$
{2k-\ell \choose \ell +1}\leq 8\cdot 2^{(1+\frac{\sqrt{5}}{5})k}
$$
by Lemma \ref{lem_binom}, hence
$$
|\cF|\leq (k^2-k+2)8^{(s-2)}2^{(1+\frac{\sqrt{5}}{5})k(s-2)}\cdot 8\cdot 2^{(1+\frac{\sqrt{5}}{5})k}=
$$
$$
= (k^2-k+2)8^{(s-1)}2^{(1+\frac{\sqrt{5}}{5})k(s-1)},
$$
which was to be proved.

{\bf Proof of Theorem \ref{main2}:}

Let $F_0\in \cF$ be a fixed subset. Clearly $|F\cap F_0|\geq \ell$  for each $F\in\cF$, since $\cF$ is an $\ell$-intersecting family. 

Let $T$ be a fixed, but arbitrary subset of $F_0$ such that $|T|=\ell$. 
Consider the set system
$$
\cF(T):=\{F\in \cF:~ T\subseteq F\}. 
$$

It is easy to see  the following Proposition.

\begin{prop} \label{Funion2}
 $$
\cF=\bigcup\limits_{T\subseteq F_0, |T|=\ell} \cF(T).
$$
\end{prop}

\qed

Define the set system
$$
\cG(T):=\{F\setminus T:~ F\in \cF(T)\}.
$$
Obviously $|\cG(T)|=|\cF(T)|$ and $\cG(T)$ is a $(k-\ell)$-uniform family.

It is easy to see that $\cG(T)$ does not contain any sunflowers with $r$ petals, because $\cF$ does not contain any sunflowers with $r$ petals. Hence for each $\alpha$
$$
|\cF(T)|=|\cG(T)|\leq D(r,\alpha)(k-\ell)! \Big( \frac{(\log \log \log (k-\ell))^2}{\alpha\log \log (k-\ell)}\Big)^{k-\ell}
$$
by Theorem \ref{Kost}.

Finally Proposition \ref{Funion} implies that
$$
|\cF|\leq \sum_{T\subseteq F_0, |T|=\ell} |\cF(T)|\leq 
$$
$$
\leq D(r,\alpha){k\choose \ell} (k-\ell)! \Big( \frac{(\log \log \log (k-\ell))^2}{\alpha \log \log (k-\ell)}\Big)^{k-\ell}.
$$
\qed

{\bf Proof of Theorem \ref{main3}:}

If $\ell:= \lceil k-\frac{k}{\log k}\rceil$, then it is easy to verify that
$$
D(r){k \choose \ell}(k-\ell)! \leq D(r)4^k.
$$


\section{Remarks}

Define $f(k,r,s)$ as the least number so that if $\cF$ is an arbitrary $k$-uniform, $L$-intersecting family, where $|L|=s$, then $|\cF|> f(k,r,s)$ implies that $\cF$ contains a sunflower with $r$ petals. In Theorem \ref{main} we proved the following recursion for $f(k,r,s)$:
$$
f(k,3,s)\leq \max_{0\leq \ell \leq k-1} {2k-\ell \choose \ell+1} f(k-1,3,s-1).
$$

Our upper bound in Theorem \ref{main} was a clear consequence of this recursion. It would be very interesting to give a similar recursion for $f(k,r,s)$ for $r>3$.

On the other hand, it is easy to prove the following Proposition from Theorem \ref{AHS}.
\begin{prop} \label{lower}
Let $1\leq s< k$ be integers. Then there exists a $c>0$ positive constant such that 
$$
f(k,3,s)> 2\cdot 10^{s/2-c\log s}
$$  
\end{prop}

\end{document}